\documentclass[12pt]{amsart}
\theoremstyle{plain}
\newtheorem{Thm}{Theorem}

\begin{document}

\title[]
{Blow-up and global solutions to $L^p$ norm preserving non-local
flows}

\author{Li Ma, Liang Cheng}

\address{Department of mathematical sciences \\
Tsinghua university \\
Beijing 100084 \\
China} \email{lma@math.tsinghua.edu.cn} \dedicatory{}
\date{Oct. 25th, 2009}

\begin{abstract}

In this paper, we study global existence and blow up properties to
$L^p$ norm preserving non-local heat flows. We first study two
kinds of $L^p$ norm preserving non-local flows and prove that
these flows have the global solutions. Finally, we give a example
to show that one kind of this heat flow may blow up in
$L^{\infty}$ norm though its $L^p$ norm is preserved.


\end{abstract}

\thanks{$^*$ The research is partially supported by the National Natural Science
Foundation of China 10631020 and SRFDP 20060003002. }
 \maketitle

\section{Introduction}
In this paper, we study global existence and blow up properties of
positive solutions to $L^p$ norm preserving non-local heat flows
$$
\partial_t u^r=\Delta
u+\lambda(t)u^s, \; \; M\times (0,T)
$$
on the Riemannian manifold $(M,g)$ with the Cauchy data, where
$T>0$, $r>0,s>0$, and $\lambda (t)$ is chosen to make the $L^p$
norm of the solution $u$ be constant. We shall show that when
$r=s=p-1>0$, the global smooth solution exists. Assume for
example, $M=\Omega$ is a bounded convex domain in $R^n$. When
$r=1$, $1<s=p<\frac{n+2}{n-2}$, the global existence of positive
is also true, however, when $r=1$ and $s=p\geq\frac{n+2}{n-2}$, we
have blowup result. Our work is motivated by the recent excellent
work of C.Caffarelli and F.Lin \cite{CL09}, where they have
studied the global existence and regularity of $L^2$ norm
preserving heat flow on bounded domains $\Omega\subset R^n$ such
as
$$
\partial_t u=\Delta
u+\lambda(t)u,
$$ with
$$
\lambda(t)=\frac{\int_{\Omega}|\nabla u|^2dx}{\int_{\Omega}u^2dx}.
$$
They also extend the method to study a family of singularly
perturbed systems of non-local parabolic equations and study the
partition problem for eigenvalues. After that the authors studied
the global existence, asymptotic behavior, stability and gradient
estimates for two kinds of non-homogenous $L^2$ norm preserving
heat flows in \cite{MC2}.

We remark that the non-local heat flow also naturally arises in
geometry such that the flow preserves some $L^p$ norm in the sense
that some the geometrical quantity (such as length, area and so
on) is preserved in the geometric heat flows. For more references
on geometric flows such as harmonic map heat flows and non-local
curve shortening flows, one may see \cite{A98}, \cite{MA},
\cite{MC} and \cite{S}.

We first study the following Yamabe type heat flow on a closed
smooth Riemannian manifold $M^n$
\begin{equation}\label{yamabe}
\left\{
\begin{array}{ll}
         u^{p-2}\partial_t u=\Delta u+\lambda(t)u^{p-1} \quad &\text{in}\ M\times\mathbb{R}_{+}, \\
          u(x,0)=g(x) &\text{in}\ M,
\end{array}
\right.
\end{equation}
where $p>1$, which has the positive solution and preserves the
$L^p$ the norm. We call equation (\ref{yamabe}) Yamabe type heat
flow since it relates to following Yamabe flow on closed manifolds
$M^n$ which introduced by Hamilton
\begin{equation}\label{hamilton_yamabe}
    \frac{\partial g}{\partial t}=(s-R)g,
\end{equation}
where $R$ denotes the scalar curvature of metric $g$ and $s$
denotes the average scalar curvature. Write
$g=u^{\frac{4}{n-2}}g_0$, $n\geq 3$, with $u$ is a positive
function and change time by a constant scale. Then
(\ref{hamilton_yamabe}) is equivalent to the following heat
equation
$$
\frac{\partial u^p}{\partial t}=L_{g_0}u+c(n)s u^p,
$$
where $p=\frac{n+2}{n-2}$ , $c(n)=\frac{n-2}{4(n-1)}$ and
$L_{g_0}u=\Delta_{g_0}u-c(n)R_{g_0}u$. For more references about
Yamabe problem and Yamabe flow, one may see \cite{B91},
\cite{B92}, \cite{H88}, \cite{Hebey}, and \cite{Y}. Now direct
computation to (\ref{yamabe}) shows that
$$
\frac{1}{p}\frac{d}{dt}\int_{M}u^pdx=\int_{M}u^{p-1} u_tdx
=-\int_{M}|\nabla u|^2dx+\lambda(t)\int_{M}u^pdx.
$$
Thus, one must have $\lambda(t)=\frac{\int_{M}|\nabla
u|^2dx}{\int_{M}g^pdx}$ to preserve the $L^p$ norm. Without loss
of generality we assume $\int_{M}g^pdx=1$. Then we consider the
following problem on closed smooth Riemannian manifold $M^n$
\begin{equation} \label{eq1}
\left\{
\begin{array}{ll}
        u^{p-2}\partial_t u=\Delta u+\lambda(t)u^{p-1} \quad &\text{in}\ M\times\mathbb{R}_{+}, \\
          u(x,0)=g(x) &\text{in}\ M,
\end{array}
\right.
\end{equation}
where $\lambda(t)=\int_{M}|\nabla u|^2dx$, $p>1$, $g(x)\geq 0$,
$\int_{M}g^pdx=1$ and $g\in C^1(M)$. Similar to the results in
\cite{Y}, we have the following theorem.

\begin{Thm}\label{thm1}
Problem (\ref{eq1}) has a positive global smooth solution $u(t)\in
L^{\infty}(\mathbb{R}_{+},H^1(M))\cap
L^{\infty}(\mathbb{R}_{+},L^{\infty}(M))$. Furthermore,
$\lambda(t)$ is non-increasing function such that $\lambda(t)\to
0$ at exponential rate as $t\to\infty$ and $u(t)$ converges
(passing by a subsequence) smoothly to a positive constant.
\end{Thm}

We next study the non-local heat flow on bounded smooth domain in
$\mathbb{R}^n$ which relates to the semilinear heat equations,
\begin{equation*}
\left\{
\begin{array}{ll}
         \partial_t u=\Delta u+\lambda(t)u^p \quad &\text{in}\ \Omega\times\mathbb{R}_{+}, \\
          u(x,0)=g(x) &\text{in}\ \Omega,\\
          u(x,t)=0 &\text{on}\ \partial\Omega\times\mathbb{R}_{+},
\end{array}
\right.
\end{equation*}
where $1<p<\frac{n+2}{n-2}$, which has the positive solution and
preserves the $L^2$ the norm. Likewise,
$$
\frac{1}{2}\frac{d}{dt}\int_{\Omega}u^2dx=\int_{\Omega}u
u_t=-\int_{\Omega}|\nabla u|^2dx+\lambda(t)\int_{\Omega}u^{p+1}dx.
$$ Thus,
one must have $\lambda(t)=\frac{\int_{\Omega}|\nabla
u|^2dx}{\int_{\Omega}u^{p+1}dx}$ to preserve the $L^2$ norm. Then
we consider the following problem on bounded smooth domain in
$\mathbb{R}^n$
\begin{equation} \label{eq2}
\left\{
\begin{array}{ll}
         \partial_t u=\Delta u+\lambda(t)u^p \quad &\text{in}\ \Omega\times\mathbb{R}_{+}, \\
          u(x,0)=g(x) &\text{in}\ \Omega,\\
          u(x,t)=0 &\text{on}\ \partial\Omega\times\mathbb{R}_{+},
\end{array}
\right.
\end{equation}
where $1<p<\frac{n+2}{n-2}$,
$\lambda(t)=\frac{\int_{\Omega}|\nabla
u|^2dx}{\int_{\Omega}u^{p+1}dx}$, $g(x)\geq 0\ \text{in}\ \Omega$,
$\int_{\Omega}g^2dx=1$ and $g\in C^1(\Omega)$. Similar to theorem
\ref{thm1}, we also have the global solution to problem
(\ref{eq2}).

\begin{Thm}\label{thm2}
Problem (\ref{eq2}) has a global positive smooth solution
$$u(t)\in L^{\infty}(\mathbb{R}_{+},H^1_0(\Omega))\cap
L^{\infty}(\mathbb{R}_{+},L^{\frac{2n}{n-2}}(\Omega)).$$ Moreover,
one can take $t_i\to \infty$ such that
$\lambda(t_i)\to\lambda_{\infty}>0$, $u(x,t_i)\to u_{\infty}(x)$
in $L^2(\Omega)$, $u(x,t_i)\rightharpoonup u_{\infty}(x)$ in
$H^1_0(\Omega)$ and $u_{\infty}$ solves the equation $\Delta
u_{\infty}+\lambda_{\infty}u_{\infty}^p=0$ in $\Omega$ and
$u_{\infty}=0$ on $\partial \Omega$ with
$\int_{\Omega}|u_{\infty}|^2dx=1$.
\end{Thm}

Finally, we find a interesting phenomenon that not all the $L^p$
norm preserving non-local flow has such good properties as problem
(\ref{eq1}). And we shall give an example following to show that
some $L^p$ norm preserving  non-local flow must blow up in
$L^{\infty}$ norm. We study the following nonlinear heat flow on
bounded smooth domain in $\mathbb{R}^n$,
\begin{equation*}
\left\{
\begin{array}{ll}
         \partial_t u=\Delta u+\lambda(t)u^p \quad &\text{in}\ \Omega\times\mathbb{R}_{+}, \\
          u(x,0)=g(x) &\text{in}\ \Omega,\\
          u(x,t)=0 &\text{on}\ \partial\Omega\times\mathbb{R}_{+},
\end{array}
\right.
\end{equation*}
where $p\geq \frac{n+2}{n-2}$, which has the positive solution and
preserves the $L^{p+1}$ the norm. Likewise,
$$
\frac{1}{p+1}\frac{d}{dt}\int_{\Omega}u^{p+1}dx=\int_{\Omega}u^p
u_t=-p\int_{\Omega}u^{p-1}|\nabla
u|^2dx+\lambda(t)\int_{\Omega}u^{2p}dx.
$$ Thus,
one must have $\lambda(t)=\frac{p\int_{\Omega}u^{p-1}|\nabla
u|^2dx}{\int_{\Omega}u^{2p}dx}$ to preserve the $L^{p+1}$ norm.
Then we consider the following problem on bounded smooth domain in
$\mathbb{R}^n$
\begin{equation} \label{eq3}
\left\{
\begin{array}{ll}
         \partial_t u=\Delta u+\lambda(t)u^p \quad &\text{in}\ \Omega\times\mathbb{R}_{+}, \\
          u(x,0)=g(x) &\text{in}\ \Omega,\\
          u(x,t)=0 &\text{on}\ \partial\Omega\times\mathbb{R}_{+},
\end{array}
\right.
\end{equation}
where $p\geq \frac{n+2}{n-2}$,
$\lambda(t)=\frac{p\int_{\Omega}u^{p-1}|\nabla
u|^2dx}{\int_{\Omega}u^{2p}dx}$, $g(x)\geq 0\ \text{in}\ \Omega$
and $g\in C^1(\Omega)$. We have the following blow up property for
problem (\ref{eq3}).

\begin{Thm}\label{thm3}
Suppose that $\Omega$ is a bounded smooth star-shaped domain in
$R^n$. Then the $L^{p+1}$ norm preserving flow (\ref{eq3}) must
blow up with $L^{\infty}$ norm in time interval $[0,\infty)$.
\end{Thm}

This paper is organized as follows. In section \ref{sect2} we
prove Theorem \ref{thm1} and Theorem \ref{thm2}. In section
\ref{sect3} we prove the blowup result, Theorem \ref{thm3}.

\section{global solutions }\label{sect2}
In this section we study the global existence property for the
$L^p$ energy preserving non-local flows. First we give the proof
of theorem \ref{thm1}.

\textbf{Proof of theorem \ref{thm1}}. Firstly by the maximum
principle, we know that $u(t) > 0$. Since
\begin{align}\label{eq1.1}
 \frac{1}{2}\frac{d}{dt}\int_{ M }|\nabla u|^2dx
 &=-\int_{ M }u_t\Delta udx\\
 &=-\int_{ M }u_t(u^{p-2}u_t-\lambda(t)u^{p-1})dx \nonumber\\
 &=-\int_{ M }u^{p-2}(u_t)^2dx+\frac{\lambda(t)}{p}\frac{d}{dt}\int_{ M }u^{p}dx
 \nonumber\\
 &=-\int_{ M }u^{p-2}(u_t)^2dx \leq 0,\nonumber
\end{align}
we know that $\lambda(t)$ is non-increasing and uniformly bounded.
we also have
\begin{align}\label{eq1.2}
 \frac{1}{2}\frac{d}{dt}\int_{ M }|\nabla u|^2dx
 &=-\int_{ M }u_t\Delta udx\\
 &=-\int_{ M }(u^{2-p}\Delta u+\lambda(t)u)\Delta u dx \nonumber\\
  &=-\int_{ M }u^{2-p}(\Delta u)^2dx+\lambda(t)\int_{ M }|\nabla u|^2 dx \nonumber
\end{align}
Note that at the maximum point of $u$, by setting
$u_{\max}(t)=\max\limits_{x\in M }(x,t)$, we have
$$
(u_{\max})_t\leq \lambda(t) u_{\max}(t).
$$
Hence
\begin{equation}\label{eq1.3}
\log\frac{u_{\max}(t)}{u_{\max}(0)}\leq \int^t_0\lambda(t)dt.
\end{equation}
Likewise, setting $u_{\min}(t)=\min\limits_{x\in M }u(x,t)$, we
have
$$
(u_{\min})_t\geq \lambda(t) u_{\min}(t).
$$
Hence
\begin{equation}\label{eq1.4}
\log\frac{u_{\min}(t)}{u_{\min}(0)}\geq \int^t_0\lambda(t)dt.
\end{equation}
Combining with (\ref{eq1.3}) and (\ref{eq1.4}), we conclude the
Harnack inequality
\begin{equation}\label{eq1.5}
u_{\max}(t)\leq C u_{\min}(t).
\end{equation}
 Since $\int_{ M }u^p(t)dx\equiv 1$, we get
\begin{equation}\label{eq1.6}
0<C'\leq u(x,t)\leq C.
\end{equation}
Now we have
\begin{equation}\label{eq1.7}
\int^{\infty}_0\lambda(t)dt<C,
\end{equation}
by (\ref{eq1.4}) and (\ref{eq1.6}). Note that the solution
$u(x,t)$ is smooth for $t>0$ by standard bootstrap argument. Hence
$\lambda(t)$ of course is continuous and problem (\ref{eq1}) has a
global solution. Now one can take a sequence $\lambda(t_i)$ such
that $\lambda(t_i)\to 0$ as $t_i\to \infty$. Moreover, since
$\lambda(t)$ is non-increasing, we know that $\lambda(t)\to 0$ as
$t\to\infty$. Furthermore, by (\ref{eq1.2}), (\ref{eq1.6}) and the
poincare inequality, we conclude that
\begin{equation}\label{eq1.8}
\lambda(t)\leq \lambda(0)\exp(-Ct).
\end{equation}
Now we integrate (\ref{eq1.1}) with t, we get
$$
\int^{\infty}_0\int_{ M }u^{p-2}(u_t)^2dx\leq C.
$$
Hence we can take a subsequence $\{t_i\}$ with $t_i\to \infty$
such that $u_i(x)=u(x,t_i)$ and we have
\begin{equation*}
\left\{
\begin{array}{ll}
         u_i\to u_{\infty}\quad &\text{in}\ L^{p}(M), \\
         u_i\rightharpoonup u_{\infty} & \text{in}\ H^1_0(M),\\
         \partial_t u_i\to 0 & \text{in}\ L^2(M).\\
\end{array}
\right.
\end{equation*}
Note that $u_{\infty}\in H^1_0(M)$ solves the equation $\Delta
u_{\infty}=0$ in $M$ and satisfies $\int_{M}|u_{\infty}|^{p}dx=1$.
Hence $u_{\infty}$ must be a positive constant. Combine with
(\ref{eq1.5}), (\ref{eq1.6}) and (\ref{eq1.8}), one can use the
same argument in \cite{Y} theorem 1 to prove $u(t_i)$ converges to
$u_{\infty}$ in $C^{\infty}$ sense. We omit the details here.
 $\Box$

Next we give the proof of theorem \ref{thm2}.

\textbf{Proof of theorem \ref{thm2}.} Firstly by the maximum
principle, we know that $u(t) > 0$. Note that
\begin{align}\label{eq2.1}
 \frac{1}{2}\frac{d}{dt}\int_{\Omega}|\nabla u|^2dx
 &=-\int_{\Omega}u_t\Delta udx\\
 &=-\int_{\Omega}u_t(u_t-\lambda(t)u^{p})dx \nonumber\\
 &=-\int_{\Omega}(u_t)^2dx+\frac{\lambda(t)}{p+1}\frac{d}{dt}\int_{\Omega}u^{p+1}dx.
 \nonumber
\end{align}
we also have
\begin{align}\label{eq2.2}
 \frac{1}{2}\frac{d}{dt}\int_{\Omega}|\nabla u|^2dx
 &=-\int_{\Omega}u_t\Delta udx\\
 &=-\int_{\Omega}(\Delta u+\lambda(t)u^p)\Delta u dx \nonumber\\
  &=-\int_{\Omega}(\Delta u)^2+p\lambda(t)\int_{\Omega}u^{p-1}|\nabla u|^2 dx. \nonumber
\end{align}
We denote that $B=\int_{\Omega}u^{p+1}dx$. Then we have
$$
\int_{\Omega}(u_t)^2dx+\frac{1}{2}\frac{d}{dt}(\lambda
B)=\frac{\lambda}{p+1}\frac{d}{dt}B,
$$
hence
$$
\frac{2\int_{\Omega}(u_t)^2dx}{\lambda
B}+\frac{d}{dt}(\log(\lambda B^{\frac{p-1}{p+1}}))=0.
$$
This implies that $\lambda B^{\frac{p-1}{p+1}}$ is non-increasing
and hence
\begin{equation}\label{eq2.3}
\lambda B^{\frac{p-1}{p+1}}(t)\leq C.
\end{equation}
 By H\"{o}lder inequality, we have
 \begin{equation}\label{eq2.4}
B(t)\geq c_0(\int_{\Omega}u^2dx)^{\frac{p+1}{2}}=c_0.
\end{equation}
By (\ref{eq2.3}), we have
\begin{equation}\label{eq2.5}
    \lambda(t)\leq C.
\end{equation}
Furthermore, by (\ref{eq2.3}), we conclude that
\begin{equation}\label{eq2.6}
    ||\nabla u||_2\leq C||u||_{p+1}.
\end{equation}
Note that $1<p<\frac{n+2}{n-2}$. Hence by Sobolev inequality, we
get
\begin{equation}\label{eq2.7}
    ||u||_{\frac{2n}{n-2}}\leq C||\nabla u||_2\leq C||u||_{p+1}\leq
    C ||u||^{\theta}_2
    ||u||^{1-\theta}_{\frac{2n}{n-2}}=C||u||^{1-\theta}_{\frac{2n}{n-2}},
\end{equation}
where $\theta=\frac{n}{(p+1)(n-1)}$. So we have
\begin{equation}\label{eq2.8}
||u||_{p+1}\leq C ||u||_{\frac{2n}{n-2}} \leq C.
\end{equation}
We integrate (\ref{eq2.1}) with t, we get
$$
\int^{\infty}_0\int_{ M }(u_t)^2dx\leq C.
$$
Note that the solution $u(x,t)$ is smooth for $t>0$ by standard
bootstrap argument. Hence $\lambda(t)$ of course is continuous and
problem (\ref{eq2}) has a global solution. Now we can take a
subsequence $\{t_i\}$ with $t_i\to \infty$ such that
$u_i(x)=u(x,t_i)$ and we have
\begin{equation*}
\left\{
\begin{array}{ll}
         u_i\to u_{\infty}\quad &\text{in}\ L^{2}(\Omega), \\
         u_i\rightharpoonup u_{\infty} & \text{in}\ H^1_0(\Omega),\\
         \partial_t u_i\to 0 & \text{in}\ L^2(\Omega),
\end{array}
\right.
\end{equation*}
by (\ref{eq2.5}), (\ref{eq2.6}) and (\ref{eq2.8}). Hence one can
take $t_i\to \infty$ such that $\lambda(t_i)\to\lambda_{\infty}$,
$u(x,t_i)\to u_{\infty}(x)$ in $L^2(\Omega)$,
$u(x,t_i)\rightharpoonup u_{\infty}(x)$ in $H^1_0(\Omega)$ and
$u_{\infty}$ solves the equation $\Delta
u_{\infty}+\lambda_{\infty}u_{\infty}^p=0$ in $\Omega$ and
$u_{\infty}=0$ on $\partial \Omega$ with
$\int_{\Omega}|u_{\infty}|^2dx=1$. Note that $\lambda_{\infty}\neq
0$ since $\Delta u_{\infty}=0$ only has zero solution in this case
which contradict to $\int_{\Omega}|u_{\infty}|^2dx=1$. $\Box$

\section{blow up}\label{sect3}
This section is devoted to the proof of theorem \ref{thm3}. Our
proof is based on the observation that the flow (\ref{eq3}) would
converge to elliptic equation $\Delta u+\lambda u^p=0$ in $\Omega$
and $u=0$ on $\partial \Omega$, $p\geq \frac{n+2}{n-2}$, with
$\int_{\Omega}|u|^{p+1}dx=1$ if we assume the $L^{\infty}$ norm is
uniformly bounded. But this equation only has the vanishing
solution if $\Omega$ is bounded smooth star-shaped domain.

\textbf{Proof of theorem \ref{thm3}.} Firstly by the maximum
principle, we know that $u(t) > 0$. Note that
\begin{align*}
 \frac{1}{2}\frac{d}{dt}\int_{\Omega}|\nabla u|^2dx
 &=-\int_{\Omega}u_t\Delta udx\\
 &=-\int_{\Omega}u_t(u_t-\lambda(t)u^{p})dx \nonumber\\
 &=-\int_{\Omega}(u_t)^2dx+\frac{\lambda(t)}{p+1}\frac{d}{dt}\int_{\Omega}u^{p+1}dx
 \nonumber\\
  &=-\int_{\Omega}(u_t)^2dx\leq 0.
 \nonumber
\end{align*}
Hence, we have
\begin{equation}\label{eq3.1}
  ||u||_{H^1}\leq C.
\end{equation}
 Now we argue by contradiction, supposing that
$||u||_{\infty}$ is uniformly bounded on time interval
$[0,+\infty)$. Since $\int_{\Omega}u^{2p}dx\geq C
\int_{\Omega}u^{p+1}dx=C$, we have
\begin{align}\label{eq3.2}
   0 \leq \lambda(t)=p\frac{\int_{\Omega}u^{p-1}|\nabla
    u|^2dx}{\int_{\Omega}u^{2p}dx}\leq Cp \int_{\Omega}u^{p-1}|\nabla
    u|^2dx \leq C||u||^2_{H^1}\leq C.
\end{align}
Integrate (\ref{eq3.1}) with t, we get
\begin{align}\label{eq3.3}
\int^t_0\int_{\Omega}(u_t)^2dx= \frac{1}{2}\int_{\Omega}|\nabla
u(t)|^2dx-\frac{1}{2}\int_{\Omega}|\nabla g|^2dx.
\end{align}
Hence
\begin{equation}\label{eq3.4}
    \int^{\infty}_0\int_{\Omega}(u_t)^2dx\leq C.
\end{equation}
Note that the solution $u(x,t)$ is smooth for $t>0$ by standard
bootstrap argument. Hence $\lambda(t)$ of course is continuous and
problem (\ref{eq3}) has a global solution. Then we can take a
subsequence $\{t_i\}$ with $t_i\to \infty$ such that
$u_i(x)=u(x,t_i)$, $\lambda(t_i)\to\lambda_{\infty}$. By
(\ref{eq3.1}), (\ref{eq3.4}) and the assumption $||u||_{\infty}$
is uniformly bounded on time interval $[0,+\infty)$, we have
\begin{equation*}
\left\{
\begin{array}{ll}
         u_i\to u_{\infty}\quad &\text{in}\ L^{p+1}(\Omega), \\
         u_i\rightharpoonup u_{\infty} & \text{in}\ H^1_0(\Omega),\\
         \partial_t u_i\to 0 & \text{in}\ L^2(\Omega).\\
\end{array}
\right.
\end{equation*}
Hence  $u_{\infty}\in H^1_0(\Omega)$ solves the equation $\Delta
u_{\infty}+\lambda_{\infty}u_{\infty}^p=0$ in $\Omega$ and
satisfies $\int_{\Omega}|u_{\infty}|^{p+1}dx=1$. This contradict
to the fact the equation
\begin{equation*}
\left\{
\begin{array}{ll}
         \Delta u+\lambda u^p=0 \quad &\text{in}\ \Omega, \\
          u(x)=0 &\text{on}\ \partial\Omega,
\end{array}
\right.
\end{equation*}
where $\lambda\geq 0$ and $p>\frac{n+2}{n-2}$, only has the
solution $u\equiv 0$ in bounded smooth star-shaped domain (see
\cite{S}).
 $\Box$

\end{document}